\documentclass[11pt,reqno]{amsart}
\usepackage{amsmath,amsfonts, amssymb, amsthm}
\usepackage{mathrsfs}
\usepackage{graphicx, color}
\usepackage{hyperref}
\usepackage{fancyhdr}
\usepackage{array}

\def\qed{\hfill \rule{4pt}{7pt}}
\def\qed{\hfill \rule{4pt}{7pt}}
\newcolumntype{A}{>{\centering\arraybackslash}p{1.9cm}}
\newcolumntype{B}{>{\centering\arraybackslash}p{2.5cm}}
\newcolumntype{C}{>{\centering\arraybackslash}p{2.1cm}}
\newcolumntype{D}{>{\centering\arraybackslash}p{0.7cm}}
\newcolumntype{E}{>{\centering\arraybackslash}p{0.7cm}}
\newcolumntype{F}{>{\centering\arraybackslash}p{2.3cm}}

\begin{document}

	\title[IDENTITIES INVOLVING PARTITIONS WITH DISTINCT ODD PARTS ]
{IDENTITIES INVOLVING PARTITIONS WITH DISTINCT ODD PARTS AND NO PARTS CONGRUENT TO 2 MOD 4}
	
	\author[Y.-C. Shen]{Yong-Chao Shen}
	\address{School of Mathematics and KL-AAGDM, Tianjin University, Tianjin, 300350, China}
	\email{\tt yc\_s@tju.edu.cn}

	\keywords{distinct odd parts, $q$-series, combinatorial proofs, restricted partitions
\newline \indent 2020 {\it Mathematics Subject Classification}. Primary 05A19, 11B65; Secondary 33B15.
\newline \indent Supported by National Key Research and Development Program of China (2023YFA1009401).}
	
	\begin{abstract} 
Recently, Pankaj Jyoti Mahanta and  Manjil P. Saikika proved some identities relating certain restricted partitions into distinct odd parts with the partition whose odd parts are distinct combinatorially. They asked for the $q$-series proofs. In this paper, we give $q$-series proofs of these identities. Also, the number of partitions into distinct odd parts equals to the number of the partitions with no parts congruent to 2 mod 4, so we can get some identities. We also give combinatorial proofs of these identities.
\end{abstract}
\maketitle

\section{Introduction}

	In recent decades, the study of some partially restricted partitions has been a very hot topic in the field of integer partitions. Andrews \cite{10} proved some analytic properties about the excess of the number of the distinct even parts partitions with even rank over those with odd rank. M.D.Hirschhorn and JA.Sellers \cite{9} studied some congruence properties of the partition with odd parts distinct. Ballantine and Welch \cite{4,6} proved some recurrence relations of partitons with distinct odd parts or distinct even parts combinatorially. Ballantine and Merca \cite{8} proved some recurrence relations between partitons with dictinct odd parts and 4-regular partitions. JA.Sellers \cite{13} proved some congruences for the partition with repeated odd parts or repeated even parts. The function of the partition with distinct odd parts also has been studied widely. For example, in the works of Ballantine and Merca \cite{15}, S.-P. Cui, W. X. Gu, Z. S. Ma \cite{16}, H. Fang, F. Xue, O. X. M. Yao \cite{17}.  
	
	In this paper we are interested in the following restricted sets of partitions:
	\begin{itemize}
		\item  $\text{Pod}(n) = $ set of all partitions of $n$ with distinct odd parts;
		\item  $\text{Pod}_{>2}(n) = $ subset of $\text{Pod}(n)$ where all parts are greater than $2$;
		\item  $\text{O1}(n) = $ set of partitions of $n$ in which no odd part is repeated and the largest part is even;
		\item  $\text{O2}(n) = $ subset of $\text{O1}(n)$ where the largest part appears at least twice;
        \item $\text{O3}(n) = $ subset of $\text{O1}(n)$ where the largest part appears exactly once;
		\item  $\text{C}(n)=$ set of partitons of $n$ in which no parts congruent to 2 mod 4;
	  \end{itemize}
we denote the cardinalities of the above sets by $pod(n)$, \text{${pod}_{>2}(n)$, $ o1(n)$,}
$o2(n)$, $o3(n)$ (\cite{3} represents by DO1($n$),  DO2($n$), DO3($n$)) and $c(n)$ respectively.
	


	Recently, Andrews and El Bachraoui \cite{2} proved three interesting identities about partition into distinct even parts and 4-regular partition by using $q$- series. In a short while, Pankaj Jyoti Mahanta and  Manjil P. Saikika \cite{3} proved those identities combinatorially. And they proved three identities regarding the partitions into distinct odd parts and the partitions in which the largest part is even, the odd parts are all distinct, and the number of occurrences of the largest part has certain restriction by analogy combinatorially. Their main results are as follows.
	
	For all $n > 1$, \hypertarget{formula:1.1}
	\[
	o1(n) + o1(n - 1) = {pod}(n), \tag{1.1}
\]	
	 and for all $ n > 4$, \hypertarget{formula:1.2}
	\[
	o2(n) + o2(n - 3) = {pod}_{>2}(n), \tag{1.2}
\]	
	and for all $ n > 2$,\hypertarget{formula:1.3} 
	\[
	o3(n+2) + o3(n - 1) = {pod}(n).\tag{1.3}
\]	



In the first part of this paper, we will  use $q$-series to prove these identities. For the first two identities, we can directly apply $q$-binomial theorem to obtain the proof. For the third identity, we modify Andrews and El Bachraoui's  proof of the Theorem 3 \cite{2} and make use of certain conclusions that they have already reached to obtain the proof. 

Moreover, we remark that,
\[
\sum_{n = 0}^{\infty} {c}(n) q^n=\sum_{n = 0}^{\infty} {pod}(n) q^n = \frac{(-q;q^2)_\infty}{(q^2;q^2)_\infty} = \frac{(q^2;q^4)}{(q;q)_\infty}, \tag{1.4} 
\]  
where the $q$-shifted factorial \cite{1} is defined by \hypertarget{formula:1.5}\[
(a;q)_0 = 1, \quad (a;q)_n = \prod_{j = 0}^{n - 1} (1 - aq^j), \quad (a;q)_\infty = \prod_{j = 0}^{\infty} (1 - aq^j),    |q|<1, \tag{1.5}
\]
 from the \hyperlink{formula:1.4}{(1.4)}, we can easily obtain the following two identities.

    For all $n > 1$,\hypertarget{formula:1.6}
    \[ o1(n) + o1(n - 1) = c(n), \tag{1.6}
    \]
    and for all $ n > 2$,\hypertarget{formula:1.7}
     \[
     o3(n+2) + o3(n - 1) = c(n).\tag{1.7}
    \] 
    
    In the second part of this paper, we  provide the combinatorial proofs for these two identities. For these two identities, since O1($n$) requires that the largest part is even, O3($n$) requires that the largest part is even and its multiplicity is exactly once, then for any partition $\lambda$ $\in$ C($n$), we will discuss it according to the parity and multiplicity of the largest part and the largest repeated odd part. Conversely, for any partition $\lambda$ $\in$ O1($n$), we just need to divide the parts that are congruent to 2 mod 4  into two equal parts. For any partition $\lambda$ $\in$ O1($n-1$), add 1 to the largest part and then we just need to divide the parts that are congruent to 2 mod 4 into two equal parts. For any partition $\lambda$ $\in$ O3($n+2$) or O3($n-1$), it also has a point in the same way.

    \section{q-SERIES PROOF OF PANKAJ JYOTI MAHANTA AND  MANJIL P. SAIKIKA'S RESULTS}
    In this section, we will prove \hyperlink{formula:1.1}{(1.1)}, \hyperlink{formula:1.2}{(1.2)} and \hyperlink{formula:1.3}{(1.3)} by using $q$-series.
    
    
    The well known $q$-binomial theorem states that \[
    	\sum_{n = 0}^{\infty} \frac{(a;q)_n}{(q;q)_n} z^n = \frac{(az;q)_\infty}{(z;q)_\infty}. \tag{2.1}
    \]
    
    It enables us to derive \hyperlink{formula:1.1}{(1.1)}.
    
   \textbf{Theorem 2.1.}  (Theorem 1.4 \cite{3}) For all $n > 1$, we have
    \[
    o1(n) + o1(n - 1) = {pod}(n).
    \]	
    
   \textbf{Proof.} It is easy to get 
   \hypertarget{formula:2.2} \[\sum_{n=0}^{\infty} o1(n)q^{n} = \sum_{n=0}^{\infty} \frac{(-q;q^{2})_{n}q^{2n}}{(q^{2};q^{2})_{n}}= \frac{(-q^{3}; q^{2})_{\infty}}{(q^{2}; q^{2})_{\infty}}. \tag{2.2}
    \].
    
Then multiplying both sides by $1 + q$, we can get

\[(1+q)\sum_{n=0}^{\infty} o1(n)q^{n} = (1+q)\sum_{n=0}^{\infty} \frac{(-q;q^{2})_{n}q^{2n}}{(q^{2};q^{2})_{n}}= \frac{(-q; q^{2})_{\infty}}{(q^{2}; q^{2})_{\infty}}=\frac{(q^{2}; q^{4})_{\infty}}{(q; q)_{\infty}}.\]

Thus, we have completed the proof of Theorem 2.1.\qed

We can also make use of \hyperlink{formula:2.1}{(2.1)} to prove \hyperlink{formula:1.2}{(1.2)}.

\textbf{Theorem 2.2.} (Theorem 1.5 \cite{3}) For all $ n > 4$, we have
\[
o2(n) + o2(n - 3) = {pod}_{>2}(n).
\]	

\textbf{Proof.} It is also easy to get 
 \[\sum_{n=0}^{\infty} o2(n)q^{n} = \sum_{n=0}^{\infty} \frac{(-q;q^{2})_{n}q^{4n}}{(q^{2};q^{2})_{n}}-1= \frac{(-q^{5}; q^{2})_{\infty}}{(q^{4}; q^{2})_{\infty}}-1, \tag{2.3}
 \] and 
 \[\sum_{n = 0}^{\infty} {pod}_{>2}(n) q^n = \frac{(-q^3;q^2)_{\infty}}{(q^4;q^2)_{\infty}}. \tag{2.4}
 \]
 
 Then  multiplying both sides of \hyperlink{formula:2.3}{(2.3)} by $1 + q^{3}$, we discover
 
 \[
 (1+q^3)\sum_{n=0}^{\infty} o2(n)q^{n} =(1+q^3) \sum_{n=0}^{\infty} \frac{(-q;q^{2})_{n}q^{4n}}{(q^{2};q^{2})_{n}}-1-q^3= \frac{(-q^{3}; q^{2})_{\infty}}{(q^{4}; q^{2})_{\infty}}-1-q^3.
 \]
 
 Thus, we have completed the proof of Theorem 2.2.\qed

We can use the Andrews and Bachraoui's results \cite{2} to prove \hyperlink{formula:1.3}{(1.3)}.

\textbf{Theorem 2.3.} (Theorem 1.6 \cite{3})For all $ n > 2$, we have
\[
o3(n+2) + o3(n - 1) = {pod}(n).
\]	

\textbf{Proof.} It is easy to see 

\[
\sum_{n=0}^{\infty} o3(n)q^{n} = \sum_{n=1}^{\infty} \frac{(-q;q^{2})_{n}q^{2n}}{(q^{2};q^{2})_{n-1}}+1. \tag{2.5}
\]

We begin at

\[\sum_{n \geq 0} \frac{(q;q)_{2n} q^{4n}}{(q^{4};q^{4})_n} = (q;q)_{\infty} \sum_{n \geq 0} \frac{q^{4n}}{(q^{4};q^{4})_n (q^{2n+1};q)_{\infty}} \]

\[= (q;q)_{\infty} \sum_{n \geq 0} \frac{q^{4n}}{(q^{4};q^{4})_n} \sum_{m \geq 0} \frac{q^{2nm + m}}{(q;q)_m} \\\]

\[= (q;q)_{\infty} \sum_{m \geq 0} \frac{q^{m}}{(q;q)_m (q^{2m + 4};q^{4})_{\infty}} \\\]

\[= (q;q)_{\infty} \sum_{m \geq 0} \frac{q^{2m}}{(q;q)_{2m} (q^{4m + 4};q^{4})_{\infty}} + (q;q)_{\infty} \sum_{m \geq 0} \frac{q^{2m + 1}}{(q;q)_{2m + 1} (q^{4m + 6};q^{4})_{\infty}} \\\]

\[= \frac{(q;q)_{\infty}}{(q^{4};q^{4})_{\infty}} \sum_{m \geq 0} \frac{(q^{4};q^{4})_m q^{2m}}{(q;q)_{2m}} + \frac{(q;q)_{\infty}}{(q^{6};q^{4})_{\infty}} \sum_{m \geq 0} \frac{(q^{6};q^{4})_m q^{2m + 1}}{(q;q)_{2m + 1}} \\\]

\[= \frac{(q;q)_{\infty}}{(q^{4};q^{4})_{\infty}} \sum_{m \geq 0} \frac{(-q^{2};q^{2})_m q^{2m}}{(q;q^{2})_m} + \frac{q(q;q)_{\infty}}{(1 - q^{2})(q^{6};q^{4})_{\infty}} \sum_{m \geq 0} \frac{(-q;q^{2})_{m+1} q^{2m}}{(q^{2};q^{2})_m}. \\\]

Then we denote $\Lambda$=$\sum_{m \geq 0} \frac{(-q;q^{2})_{m + 1} q^{2m}}{(q^{2};q^{2})_m}$.

Andrews and Bachraoui \cite{2} proved \[ \sum_{m \geq 0} \frac{(-q^{2};q^{2})_m q^{2m}}{(q;q^{2})_m}=\frac{q(q^{4};q^{4})_{\infty}}{(1 + q^{3})(q;q)_{\infty}} + \frac{1 - q}{1 + q^{3}}.\]

Then we  use their results to obtain

\[\sum_{n \geq 0} \frac{(q;q)_{2n} q^{4n}}{(q^{4};q^{4})_n}=\frac{(q;q)_{\infty}}{(q^{4};q^{4})_{\infty}}\left(\frac{q(q^{4};q^{4})_{\infty}}{(1 + q^{3})(q;q)_{\infty}} + \frac{1 - q}{1 + q^{3}}\right)+\frac{q(q;q)_{\infty}}{(q^{2};q^{4})_{\infty}}\Lambda\]

\[=\frac{q}{1 + q^{3}}+\frac{(q;q)_{\infty}}{(q^{4};q^{4})_{\infty}}\frac{1 - q}{1 + q^3}+\frac{q(q;q)_{\infty}}{(q^{2};q^{4})_{\infty}}\Lambda.\]

Then multiplying both sides by $\frac{(q^{2};q^{4})_{\infty}}{q(q;q)_{\infty}}$, we get

\[\frac{(q^{2};q^{4})_{\infty}}{q(q;q)_{\infty}}\sum_{n \geq 0} \frac{(q;q)_{2n} q^{4n}}{(q^{4};q^{4})_n}\]

\[=\Lambda+\frac{(q^{2};q^{4})_{\infty}}{(1 + q^{3})(q;q)_{\infty}}+\frac{(1 - q)(q^{2};q^{4})_{\infty}}{q(1 + q^{3})(q^{4};q^{4})_{\infty}}.\]

Andrews and Bachraoui \cite{2} also proved
\[
\frac{(q^{4};q^{4})_{\infty}}{(q;q)_{\infty}}\sum_{n \geq 0} \frac{(q;q)_{2n} q^{4n}}{(q^{4};q^{4})_n}=\frac{2q(q^{4};q^{4})_{\infty}}{(1+q^{3})(q;q)_{\infty}}+\frac{1-q}{1+q^{3}}.
\]

So, we have

\[\Lambda=\frac{(q^{2};q^{4})_{\infty}}{q(q;q)_{\infty}}\sum_{n \geq 0} \frac{(q;q)_{2n} q^{4n}}{(q^{4};q^{4})_n}-\frac{(q^{2};q^{4})_{\infty}}{(1 + q^{3})(q;q)_{\infty}}-\frac{(1 - q)(q^{2};q^{4})_{\infty}}{q(1 + q^{3})(q^{4};q^{4})_{\infty}}\\\]

\[=\frac{(q^{2};q^{4})_{\infty}}{q}\left(\frac{2q}{(1 + q^{3})(q;q)_{\infty}}+\frac{1 - q}{(1 + q^{3})(q^{4};q^{4})_{\infty}}\right)-\frac{(q^{2};q^{4})_{\infty}}{(1 + q^{3})(q;q)_{\infty}}-\frac{(1 - q)(q^{2};q^{4})_{\infty}}{q(1 + q^{3})(q^{4};q^{4})_{\infty}}\\\]

\[=\frac{2(q^{2};q^{4})_{\infty}}{(1 + q^{3})(q;q)_{\infty}}-\frac{(q^{2};q^{4})_{\infty}}{(1 + q^{3})(q;q)_{\infty}}=\frac{(q^{2};q^{4})_{\infty}}{(1 + q^{3})(q;q)_{\infty}}.\\\]

Then we have
\[\sum_{n = 0}^{\infty} o3(n)q^{n}=q^{2}\Lambda + 1=\frac{q^{2}(q^{2};q^{4})_{\infty}}{(1 + q^{3})(q;q)_{\infty}}+1.\\ \tag{2.6}\]
Then multiplying both sides by $1+q^{3}$, we get

\[(1 + q^{3})\sum_{n = 0}^{\infty} o3(n)q^{n}=q^{2}\frac{(q^{2};q^{4})_{\infty}}{(q;q)_{\infty}}+1 + q^{3}.\]

It reveals that for $k>3$, $o3(k$)+$o3(k-3$)=$c(k-2$). Let $n$=$k-2$, then we have for $n>1$,

 \[o3(n+2)+o3(n-1)=c(n).\]   \qed
 
 \section{COMBINATORIAL PROOF OF OUR RESULTS}         
Recall that a partition of $n$ is a non-increasing sequence $\lambda=(\lambda_1,\lambda_2,\cdots,\lambda_k)$ such that $\lambda_1 \geq \lambda_2 \geq \cdots\geq \lambda_k $, and $\sum_{i = 1}^{k} \lambda_{i} = n$, \text{the  size  of }$\lambda$, denoted $|\lambda|$, is the sum of the parts. $m(\lambda_i)$ denotes the multiplicity of $\lambda_i$ in the partition $\lambda$. Also, we sometimes use exponent notation to represent the multiplicity of parts, i.e., $({5}^{2}, 2)$ denotes the partition $(5, 5, 2)$.

In this section, we will use some operations to construct bijections to prove \hyperlink{formula:1.6}{(1.6)} and \hyperlink{formula:1.7}{(1.7)}.

The  key operations which we will use are as follows.

Given two partitions $\lambda$ and $\mu$, we define $\lambda \cup \mu$ to be the multiset union of the parts of
$\lambda$ and $\mu$. For example, 
\[
\lambda=(5,4,3,2,1), \mu=(4,4,3,3,2),   \text{ then }                         \lambda \cup \mu=(5,4,4,4,3,3,3,2,2,1).
\]
For a partition set which odd parts appear with even multiplicity, defined by $A$, for a partition set which only has even parts, defined by $B$. We define a map $\Phi$: $A$ $\rightarrow$ $B$:
for a partition $\lambda$ $\in$ $A$, merging the equal odd parts two by two until the transformed partition has no  odd parts, such as, \[\lambda=(5,5,4,3,3,3,3,2), \text{ then } \Phi(\lambda)=(10,6,6,4,2).\]

\textbf{Theorem 3.1.} For all $n > 1$, we have
\[
o1(n) + o1(n - 1) = c(n).
\]	
\textbf{Proof.} For \(n > 1\), \(\lambda = (\lambda_1, \lambda_2, \dots, \lambda_k) \in \text{C}(n)\), we divide \(\lambda\) into seven cases.

\textbf{CASE 1:} \(\lambda_1\) is even. We write \(\lambda = \alpha \cup \beta\), where \(\alpha\) is a partition of distinct odd parts and \(\beta\) is a partition in which odd parts appear with even multiplicity. Then \(\alpha \cup \Phi(\beta) \in \text{O1}(n)\).

\textbf{CASE 2:}  \(\lambda_1 \) is odd and there is no repeated odd part in $\lambda$. Let \(\lambda' = (\lambda_1 - 1, \lambda_2, \dots, \lambda_k)\) $\in \text{O1}(n - 1)$.

\textbf{CASE 3:} \(\lambda_1\) is odd and \(m(\lambda_1)>1\). Then, by merging 2
copies of \(\lambda_1\) into a single part, we get the partition \(\mu\). We write \(\mu = \alpha \cup \beta\), where \(\alpha\) is a partition of distinct odd parts and \(\beta\) is a partition in which odd parts appear with even multiplicity. Then \(\alpha \cup \Phi(\beta) \in \text{O1}(n)\).

\textbf{CASE 4:}  \(\lambda_1 = 4k + 1\), $k \in N$ and \(m(\lambda_1)=1\). 
Assume that \(\lambda\) contains repeated odd parts. 
Suppose the largest repeated odd part is less than  \(\frac{\lambda_1 + 1}{2}\). Let \(\lambda' = (\lambda_1 - 1, \lambda_2, \dots, \lambda_k)\). We write \(\lambda' = \alpha \cup \beta\), where \(\alpha\) is a partition of distinct odd parts and \(\beta\) is a partition in which odd parts appear with even multiplicity. Then \(\alpha \cup \Phi(\beta) \in \text{O1}(n-1)\).

\textbf{CASE 5:} \(\lambda_1 = 4k + 1\), $k \in N$ and \(m(\lambda_1)=1\). 
Assume that \(\lambda\) contains repeated odd parts. 
Suppose the largest repeated odd part is greater than or equal to \(\frac{\lambda_1 + 1}{2}\). We write \(\lambda = \alpha \cup \beta\), where \(\alpha\) is a partition of distinct odd parts and \(\beta\) is a partition in which odd parts appear with even multiplicity. Then \(\alpha \cup \Phi(\beta) \in \text{O1}(n)\).

\textbf{CASE 6:} \(\lambda_1 = 4k + 3\), $k \in N$ and \(m(\lambda_1)=1\). 
Assume that \(\lambda\) contains repeated odd parts. 
Suppose the largest repeated odd part is greater than  \(\frac{\lambda_1 - 1}{2}\). We write \(\lambda = \alpha \cup \beta\), where \(\alpha\) is a partition of distinct odd parts and \(\beta\) is a partition in which odd parts appear with even multiplicity. Then \(\alpha \cup \Phi(\beta) \in \text{O1}(n)\).

\textbf{CASE 7:}  \(\lambda_1 = 4k + 3\), $k \in N$ and \(m(\lambda_1)=1\). 
Assume that \(\lambda\) contains repeated odd parts. 
Suppose the largest repeated odd part is less than or equal to \(\frac{\lambda_1 - 1}{2}\). Let \(\lambda' = (\lambda_1 - 1, \lambda_2, \dots, \lambda_k)\). We write \(\lambda' = \alpha \cup \beta\), where \(\alpha\) is a partition of distinct odd parts and \(\beta\) is a partition in which odd parts appear with even multiplicity. Then \(\alpha \cup \Phi(\beta) \in \text{O1}(n-1)\).

Conversely,

\textbf{CASE 1:} \(\lambda = (\lambda_1, \lambda_2, \dots, \lambda_k) \in \text{O1}(n)\). Divide each part that is congruent to 2 mod 4 in \(\lambda\) into two equal numbers to obtain \(\mu\).Then $\mu \in \text{C}(n)$.

\textbf{CASE 2:} \(\lambda = (\lambda_1, \lambda_2, \dots, \lambda_k) \in \text{O1}(n - 1)\), let \(\alpha=(\lambda_1 + 1, \lambda_2, \dots, \lambda_k)\). Divide each part that is congruent to 2 mod 4 in \(\alpha\) into two equal numbers to obtain \(\mu\). Then $\mu \in \text{C}(n)$.

Therefore,  $ o1(n) +  o1(n - 1) = {c}(n).$\qed

\textbf{Example 1:}   We give some examples as follows. 

\begin{table}[htbp]
	\centering  
\begin{tabular}{|c|c|c|c|c|} 
	\hline
	$\lambda$ & $\alpha$ & $\beta$ & $\Phi(\beta)$ & $\alpha \cup \Phi(\beta)$ \\
	\hline
	$({12}^{2},{11}^{3},8,3)$ & $({12}^{2},11,8,3)$ & $({11}^{2})$ & $(22)$ & $(22,{12}^{2},11,8,3)$ \\
	\hline
	$(13,12,{11}^{4},7,{3}^{2})$ & $(13,12,7)$ & $({11}^{4},{3}^{2})$ & $({22}^{2},6)$ & $({22}^{2},13,12,7,6)$ \\
	\hline
	$(11,{9}^{3},8,3,1)$ & $(11,9,8,3,1)$ & $({9}^{2})$ & $(18)$ & $(18,11,9,8,3,1)$ \\
	\hline

\end{tabular}
\vspace{0.5em}
\caption{Examples for the operation process.}  
\label{tab:example}  
\end{table}





\textbf{Example 2:} For $\lambda$=$(11,8,{5}^{3},4,3)$ $\in$ C(41), then we apply it in CASE7, we can get the \[\lambda'=(10,8,{5}^{3},4,3), \alpha=(10,8,5,4,3) \text{ and the } \beta=({5}^{2}), \Phi(\beta)=(10),\] so $\alpha$ $\cup$ $\Phi(\beta)$=$({10}^{2},8,5,4,3)$ $\in$ O1(40). 

Conversely, for partition \[\mu=({10}^{2},8,5,4,3) \in O1(40),\text{ let } \alpha=(11,10,8,5,4,3).\] We divide the parts that are congruent to 2 mod 4 in $\alpha$ into two equal numbers. We get the partition $(11,8,{5}^{3},4,3)$ $\in$ C(41).

\textbf{Theorem 3.2.} For all $ n > 2$,we have
\[
o3(n+2) + o3(n - 1) = c(n).
\]
\textbf{Proof.} For \(n >2\), \(\lambda = (\lambda_1, \lambda_2, \dots, \lambda_k) \in \text{C}(n)\), we divide \(\lambda\) into eleven cases.

\textbf{CASE 1:}  $\lambda_1 $ \text{is even, the odd parts are distinct}, then $\mu$= $\phi_1 (\lambda)$ $ =(\lambda_1+2, \lambda_2, \dots, \lambda_k) \in \text{O3}(n+2)$.

\textbf{CASE 2:}  $\lambda_1$ is even, $\lambda$ has the repeated odd parts, and the largest repeated odd part is greater than or equal to $\frac{\lambda_1+2}{2}$, denoted as \(a\). Merge two parts of value \(a\) and add 2 to form a single part, then obtain \(\mu\).
Then, write \(\mu=\gamma \cup \beta\), where \(\gamma\) is the partition with distinct odd parts, and \(\beta\) is a partition in which odd parts appear with even multiplicity. Then $\Phi(\beta) \cup \gamma\in \text{O3}(n + 2)\).

\textbf{CASE 3:}   $\lambda_1$ is even, $\lambda$ has the repeated odd parts, and the largest repeated odd part is less than \(\frac{\lambda_1 + 2}{2}\). Let \(\mu = (\lambda_1 + 2, \lambda_2, \dots, \lambda_k)\). We write  \(\mu=\gamma \cup\beta\), where \(\gamma\) is the partition with distinct odd parts, and \(\beta\) is a partition in which odd parts appear with even multiplicity. Then $\Phi(\beta) \cup \gamma \in \text{O3}(n + 2)\).

\textbf{CASE 4:}  $\lambda_1$ is odd and there is no repeated odd part in $\lambda$, and $\lambda_2<\lambda_1-1$, then $\mu$= $\phi_2(\lambda)=(\lambda_1-1, \lambda_2, \dots, \lambda_k)$ $\in \text{O3}(n-1)$.

\textbf{CASE 5:}  $\lambda_1$ is odd and there is no repeated odd part in $\lambda$, and $\lambda_2=\lambda_1-1$,
then $\mu=$$\phi_3(\lambda)=(\lambda_1+1, \lambda_2+1, \dots, \lambda_k)$ $\in \text{O3}(n+2).$

\textbf{CASE 6:} 
 $\lambda_1$ is odd and \(m(\lambda_1)\geq 2\). 
Merge two parts of value \(\lambda_1\) and add 2 to form a single part, then obtain \(\mu\). 
Write \(\mu=\gamma \cup \beta\), where \(\gamma\) is the partition with distinct odd parts, and \(\beta\) is a partition in which odd parts appear with even multiplicity. 
Then $\Phi(\beta) \cup  \gamma \in \text{O3}(n + 2)\).

\textbf{CASE 7:} \(\lambda_1 = 4k + 1\), $k \in N$ and \(m(\lambda_1)=1\). 
Assume that \(\lambda\) contains repeated odd parts. 
Suppose the largest repeated odd part is greater than or equal to \(\frac{\lambda_1 + 1}{2}\), 
denoted as \(a\). Merge two parts of value \(a\) and add 2 to form a single part, then obtain \(\mu\). Write \(\mu=\gamma \cup \beta\), where \(\gamma\) is the partition with distinct odd parts, and \(\beta\) is a partition in which odd parts appear with even multiplicity. Then $\Phi(\beta) \cup \gamma\in \text{O3}(n + 2)\).

\textbf{CASE 8:}  \(\lambda_1 = 4k + 1\), $k \in N$ and \(m(\lambda_1)=1\). 
Assume that \(\lambda\) contains repeated odd parts.  
Suppose the largest repeated odd part is less than \(\frac{\lambda_1 + 1}{2}\) and $m(\lambda_1 - 1)$   is equal to 0. Let $\mu = (\lambda_1 - 1, \lambda_2, \dots, \lambda_k)$, write \(\mu=\gamma \cup \beta\), where $\gamma$ is the partition with distinct odd parts, and \(\beta\) is a partition in which odd parts appear with even multiplicity.  
Then $\Phi(\beta) \cup \gamma \in \text{O3}(n -1)\).

\textbf{CASE 9:}  \(\lambda_1 = 4k + 1\), $k \in N$ and \(m(\lambda_1)=1\). 
Assume that \(\lambda\) contains repeated odd parts.  
Suppose the largest repeated odd part is less than $\frac{\lambda_1 + 1}{2}$ and $m(\lambda_1 - 1)$ $\geq$ 1, then $\lambda_1 $ add  1, select a $ \lambda_1-1 $ add  1 obtain  $\mu$.
Write $\mu=\gamma \cup \beta$, where $\gamma$ is the partition with distinct odd parts, and \(\beta\) is a partition in which odd parts appear with even multiplicity. 
Then $\Phi(\beta) \cup \gamma \in \text{O3}(n + 2)\).

\textbf{CASE 10:} 
 \(\lambda_1 = 4k + 3\), $k \in N$ and \(m(\lambda_1)=1\). 
Assume that \(\lambda\) contains repeated odd parts. 
Suppose the largest repeated odd part is greater than or equal to \(\frac{\lambda_1 - 1}{2}\), 
denoted as \(a\). Merge two parts of value \(a\) and add 2 to form a single part, then obtain \(\mu\). Write \(\mu=\gamma \cup \beta\), where \(\gamma\) is the partition with distinct odd parts, and \(\beta\) is a partition in which odd parts appear with even multiplicity. Then $\Phi(\beta) \cup \gamma\in \text{O3}(n + 2)\).

\textbf{CASE 11:} 
 \(\lambda_1 = 4k + 3\), $k \in N$ and \(m(\lambda_1)=1\). 
Assume that \(\lambda\) contains repeated odd parts.
Suppose the largest repeated odd part is  less than $\frac{\lambda_1 - 1}{2}$. Let $\mu=(\lambda_1-1, \lambda_2, \dots, \lambda_k)$, write $\mu=\gamma \cup \beta$, where $\gamma$ is the partition with distinct odd parts, and \(\beta\) is a partition in which odd parts appear with even multiplicity. 
Then $\Phi(\beta) \cup \gamma \in \text{O3}(n -1)\).

Conversely,

\textbf{CASE 1:} Let \(\lambda = (\lambda_1, \lambda_2, \dots, \lambda_k) \in \text{O3}(n - 1)\).  \(\alpha =\phi_3(\lambda)= (\lambda_1 + 1, \lambda_2, \dots, \lambda_k)\). Divide each part that is congruent to 2 mod 4 in \(\alpha\) into two equal numbers to obtain \(\mu\). Then, \(\mu \in \text{C}(n)\).

\textbf{CASE 2:} Let \(\lambda = (\lambda_1, \lambda_2, \dots, \lambda_k) \in \text{O3}(n +2)\).  \(\alpha =\phi_2(\lambda)= (\lambda_1 -2, \lambda_2, \dots, \lambda_k)\). Divide each part that is congruent to 2 mod 4 in \(\alpha\) into two equal numbers to obtain \(\mu\). Then, \(\mu \in \text{C}(n)\).

Therefore $ o3(n+2) +  o3(n - 1) = {c}(n)$. \qed

\textbf{Example 3} We give some examples as follows. 

\begin{table}[htbp]
	\centering 
\begin{tabular}{|A|B|C|D|E|F|}
	\hline
	$\lambda$ & $\mu$ & $\gamma$ & $\beta$ & $\Phi(\beta)$ & $\gamma \cup \Phi(\beta)$\\ 
	\hline
	$({12}^{2},{11}^{4},8)$ & $(24,{12}^{2},{11}^{2},8)$ & $(24,{12}^{2},8)$ & $({11}^{2})$  & $(22)$ & $(24,22,{12}^{2},8)$  \\
	\hline
	$(13,{11}^{3},{3}^{2})$ & $(24,13,11,{3}^{2})$ & $(24,13,11)$ & $({3}^{2})$ & $(6)$ & $(24,13,11,6)$  \\
	\hline
	$(11,{8}^{2},5,{3}^{3})$ & $(10,{8}^{2},5,{3}^{3})$ & $(10,{8}^{2},5,3)$ & $({3}^{2})$ & $(6)$ & $(10,{8}^{2},6,5,3)$ \\
	\hline
	$(9,{8}^{2},5,{3}^{3})$ & $(10,9,8,5,{3}^{3})$ & $(10,9,8,5,3)$ & $({3}^{2})$ & $(6)$ & $(10,9,8,6,5,3)$ \\
	\hline
	\end{tabular}
	\vspace{0.5em}
	\caption{Examples for the operation process.}  
	\label{tab:example}  
	\end{table}

    \end{document}